\documentclass[12pt,a4paper,leqno]{amsart}

\usepackage{amsmath}
\usepackage{amssymb}
\swapnumbers
\theoremstyle{plain}

\newtheorem{lemm}{Lemma}[section]
\newtheorem{pro}[lemm]{Proposition}

\theoremstyle{definition}

\theoremstyle{remark}
\newtheorem{remk}[lemm]{Remark}

\numberwithin{equation}{section}

\newcommand{\thismonth}{\ifcase\month\or
  January\or February\or March\or April\or May\or June\or July\or
  August\or September\or October\or November\or December\fi
  \space\number\year}
\DeclareMathAlphabet{\mathrmsl}{OT1}{cmr}{m}{sl}

\newcommand{\oper}[3][n]{\newcommand{#2}{\mathop{\mathrm{#3}}%
\ifx n#1\nolimits\else\limits\fi} }
\newcommand{\rsoper}[3][n]{\newcommand{#2}{\mathop{\mathrmsl{#3}}%
\ifx n#1\nolimits\else\limits\fi} }
\newcommand{\BM}{{\mathbb B}}

\newcommand{\ZM}{{\mathbb Z}}
\newcommand{\RM}{{\mathbb R}}
\newcommand{\CM}{{\mathbb C}}
\newcommand{\ra}{\rangle}
\newcommand{\la}{\langle}
\newcommand{\dn}{\mathcal{D}}
\newcommand{\dl}{\widetilde{\mathcal{D}}}
\newcommand{\bu}{\cdot}

\newcommand{\iden}{\mathrm{id}}
\newcommand{\nuv}{\widehat{\nabla}}
\renewcommand{\div}{\mathrm{div}}
\renewcommand{\ker}{\mathrm{Ker}}

\newcommand{\Scal}{\mathrm{Scal}}
\newcommand{\Ric}{\mathrm{Ric}}
\newcommand{\spin}{\mathrm{Spin}}
\newcommand{\spic}{\mathrm{Spin}^c}

\renewcommand{\geq}{\geqslant}
\renewcommand{\leq}{\leqslant}
\newcommand{\proofof}[1]{\end{#1}\begin{proof}}

\begin{document}

\title[Asymptotically complex hyperbolic spaces]{Rigidity at infinity for 
even-dimensional asymptotically complex hyperbolic spaces}
\author{Hassan Boualem}
\author{Marc Herzlich}
\address{D\'epartement des Sciences Math\'e\-matiques\\ 
G\'eom\'etrie - Topologie et
Al\-g\`ebre\\ UMR 5030 du CNRS\\ Universit\'e
Montpellier~II\\ France}
\thanks{The second author is a member of the ``EDGE'' Research Training 
Network HPRN-CT-2000-00101 supported by the Human Potential Program of the
European Union.}

\begin{abstract}
Any K\"ahler metric on the ball which is strongly asymptotic to 
complex hyperbolic space and whose scalar curvature is no less than the 
one of the complex hyperbolic space must be isometrically biholomorphic 
to it. This result has been known for some time 
in odd complex dimension and we provide here a proof in even dimension. 
\end{abstract}


\maketitle

\bigskip

\section*{Introduction.}

{\flushleft\it Motivations}. -- In \cite{mh-hyper}, the second author 
proved that any asymptotically complex-hyper\-bolic 
spin and K\"ahler manifold {\sl of odd complex dimension}
whose scalar curvature is larger
than the one of complex hyperbolic space cannot be too close at infinity 
to the model space. This stood as a K\"ahlerian counterpart of the 
analogous rigidity statements proved by L.~Andersson and M.~Dahl 
\cite{dahl2}, M.~Min-Oo \cite{minoo-hyperbolic}
and M.~C.~Leung \cite{leung-pinching} for the case of the real hyperbolic
space. It also appeared reminsicent of well-known rigidity results for
asymptotically flat spaces, known as positive mass theorems. The purpose of 
this short note is to extend the K\"ahlerian rigidity 
statement to the (up to now missing) even-dimensional (complex) case.

\medskip

{\flushleft\it Notations}. -- We consider in this paper complete K\"ahler 
manifolds $(M,g,J)$ of even complex dimension $m=2n$. Its K\"ahler
form will always be denoted by $\omega$. We will denote by $R$ its Riemann 
curvature tensor, $\Ric$ its Ricci tensor, with $\Scal$ its scalar 
curvature and $\Ric_0$ its trace-free part, and $\rho$ the Ricci $2$-form, 
$\rho_0$ being its trace-free part with respect to the K\"ahler form.

An example is the complex hyperbolic space $\CM \mathbf{H}^{2n}$ which is 
the unique simply connected K\"ahler manifold with constant negative 
holomorphic sectional curvature (it will be normalized to $-4$ 
throughout this paper). It is diffeomorphic to the unit ball $\BM$
in $\CM^n$ and
is Einstein with scalar curvature $-4m(m+1)$.

Any K\"ahler manifold $(M,g,J)$ is said to be {\sl strongly asymptotic} to 
complex hyperbolic space if  

\smallskip

\begin{enumerate}
\item[{\bf(i)}] there are a compact set $K$ of $M$ and
a diffeomorphism between 
$M\setminus K$ and the exterior of a ball in $\CM\mathbf{H}^m$, so that we 
can 
consider $M\setminus K$ to be endowed with two K\"ahler structures: the 
original $(g,J)$ and the standard $(g_0,J_0)$;
\item[{\bf (ii)}] if $A$ is defined on $M\setminus K$ by
$g(X,Y) =  g_0(AX,AY)$ for all $X$, $Y$, then $A$, $J_0J^{-1}$ and 
$J^{-1}J_0$ are uniformly bounded from below and above; 
\item[{\bf(iii)}] there is an $\varepsilon>0$ such that, if 
$r=d_{g_{0}}(x_{0},.)$ for an arbitrary basepoint $x_{0}$ in $M$,
$|\nabla^{g_{0}}A| + |A-I|
+ |J^{-1}J_{0}-I| + |J_{0}J^{-1}-I| = O(e^{-(2m+2+\varepsilon)r})$;
\item[{\bf (iv)}] the scalar curvature of $(M,g)$ is uniformly bounded.
\end{enumerate}

\medskip

{\flushleft In this paper, we prove :}

\smallskip

{\flushleft\bf Theorem}. {\sl Let} $(\BM,g,J)$ {\sl be a K\"ahler
metric on the ball of} even {\sl complex dimension} $m$.
{\sl If} $(g,J)$ {\sl is strongly asymptotic to the complex 
hyperbolic metric and} $\Scal_g \geq -4m(m+1)$,
{\sl then} $(\BM,g,J)$ {\sl is biholomorphically isometric to} 
$\CM \mathbf{H}^{m}$.

\medskip

The proof in \cite {mh-hyper} relied 
heavily on the existence on each complex hyperbolic space $\CM\mathbf{H}^m$ 
of {\sl odd} complex dimension of a special set of spinor fields, called 
{\sl K\"ahlerian Killing spinors}. The unfortunate non-existence of 
such spinor fields in {\sl even} complex dimension (which is in a sense
the non-compact counterpart of the well-known fact, in the realm of compact
manifolds, that complex projective spaces $\CM\mathbf{P}^{2n}$ are not
spin) was the cause for the dimensional restriction. 

\smallskip

The main contribution of this paper is that an analogous proof for the 
even dimensional case may be available, {\sl once the right objects are
used}. It is based on a simple (but crucial) fact, proven below, that 
seems to have escaped notice so far : whereas distinguished spinor fields 
do not exist for 
the classical $\spin$--structure on $\CM\mathbf{H}^{2n}$, there does exist 
a lot of distinguished sections of spinor bundles issued from well-chosen
$\spin^c$--structures. 
This remark paves the way for using the techniques already developed in 
\cite{mh-hyper}, but shifting from the $\spin$ to the $\spin^c$-context 
renders the proofs more involved and also unfortunately makes the technique 
slightly less adequate. 

\smallskip

Although we restricted ourselves to the ball in
the Theorem above, 
we will prove below a slightly stronger statement, which is our main 
technical result and implies immediately the previous one.

\smallskip

{\flushleft\bf Theorem} (technical version).
{\sl Let} $(M,g,J)$ {\sl be a K\"ahler
manifold of} even {\sl complex dimension} $m$, {\sl such that} 
$\left[\frac{\omega}{i\pi}\right]$ 
{\sl is an integral class and its associated
line bundle defines a} $\spic$--{\sl structure}.
{\sl If} $(M,g,J)$ {\sl is strongly asymptotic to the complex 
hyperbolic metric and} $\Scal_g \geq -4m(m+1)$,
{\sl then} $(M,g,J)$ {\sl is biholomorphically isometric to} 
$\CM \mathbf{H}^{m}$.

\medskip

{\flushleft\it Remarks}. 
1.-- Although the topological conditions on the 
K\"ahler class appear essentially as {\sl ad hoc} assumptions for the 
Theorem to hold, they are satisfied in a large number of situations. 
They are for instance preserved by blow-up.

2.-- Apart from the analogy with the odd-dimensional
case, another motivation is the construction by C.~R.~Graham and J.~Lee 
\cite{graham-lee} and more recently by O.~Biquard \cite{biquard-asympt-sym} 
of a large family of complete Einstein metrics on the ball. These are 
deformations of the standard metrics of non-compact rank-one symmetric 
spaces and are in one-to-one correspondance with a neighbourhood of 
deformations of the structure at infinity of each symmetric space. This 
favours the idea that each admissible structure at infinity can be 
filled by a unique Einstein manifold (possibly not topologically trivial). 
Graham-Lee and Biquard results settle the {\sl local} (close to the 
standard 
structure) existence and uniqueness problem, whereas results as in 
\cite{mh-hyper,leung-pinching,minoo-hyperbolic} provide {\sl global}
uniqueness, at least when the structure at infinity is standard enough.
Unfortunately, even the technical version of our result brings only very 
partial information in this context. 

\smallskip

{\flushleft\it Contents of the note}. In section 1, we describe the model
spinors on the complex hyperbolic space and begin the proof 
of the Theorem (in its technical version): we show that any manifold 
satisfying the assumptions bears K\"ahlerian Killing spinors for an 
adequate $\spin^c$--structure.
Section 2 is then devoted to the geometrical constructions leading
to the final rigidity statement: any such manifold must be the complex
hyperbolic space. Since parts of the proof are simple reproduction of 
arguments given in \cite{mh-hyper}, we have chosen to be quite short for 
these steps and to stress only the points that differ from the previously 
existing proof.

\medskip

\section{The fundamental formula and its consequences}

From now one, we will study a K\"ahler manifold $(M,g,J)$ satisfying 
all the assumptions of the theorem above. If $\frac{\omega}{i\pi}$ is in 
the image of $H^2(M,\ZM)$ in $H^2(M,\RM)$, $F = -2i\omega$ is the curvature 
form of a hermitian connection on a complex line bundle $L$. Its associated 
$S^1$-bundle will be denoted by $\pi: P \rightarrow M$. We consider the 
$\spic$--structure on $M$ induced by this choice as an auxiliary bundle.
Its spinor bundle will be denoted by $\Sigma^c$ and it splits under the 
action of the K\"ahler form as 
\[ \Sigma^c = \bigoplus_{0\leq q \leq m} \ker \left(\omega\bu + i(m-2q) 
\iden\right) = \oplus \ \Sigma_q^c.\]
A K\"ahlerian Killing spinor is a section $\Psi = \psi_{r-1} + \psi_r$ 
(for some $r$ in $\{0,...,m\}$) solving 
\[ 
\nabla_X\psi_{r-1} + i X^{0,1}\cdot\psi_r = 0,\ \   
\nabla_X\psi_r + i X^{1,0}\cdot\psi_{r-1} = 0.
\]
For future reference, we note that, equivalently 
\cite{kirchberg-killing-kahler},
\begin{equation}\label{eq:kki2}
\nabla_{X}\Psi + \frac{i}{2}X\bu\Psi + \frac{(-1)^{n}}{2}JX\bu
\overline{\Psi}\ =\ 0\ \ \ \ \forall\, X \in TM.
\end{equation} 

We begin with a quick tour of the model space. Letting $K$ be the canonical
bundle of the complex hyperbolic space, the elementary but fundamental 
remark is the following 

\begin{pro} \label{pro:exiskki}
Let $m=2n$. Then the spinor subbundles $\Sigma^c_{n-1} 
\oplus \Sigma_n^c$ of the $\spic$--structure on $\CM\mathbf{H}^m$ induced 
by the choice $L=K^{-\frac{1}{m+1}}$ 
are trivialized by K\"ahlerian Killing spinors. \end{pro}

\begin{proof} We relying on the explicit computations done in 
\cite{kirchberg-killing-kahler}. Any spinor field $\psi\in\Sigma^c_{q}$ 
can be written as $\bar\omega\otimes\psi_0$ where $\omega$ is a form of 
type $(q,0)$ on $\CM\mathbf{H}^m$ and $\psi_0$ is a holomorphic section 
of $\Sigma^c_{0}=K^{\frac{n}{m+1}}$. Then the set of equations 
\begin{equation}\label{for:explicitkki}
\nabla_{X^{0,1}} \left( |\psi_0|^2 \omega \right) =0, \ \textrm{ and }\ 
\nabla_{X^{1,0}}\omega = \frac{1}{n}\, \iota_{X^{1,0}}(\partial\omega)
\ \ \textrm{ for any } X,
\end{equation}
is equivalent to  
\begin{equation}\label{for:explicitkk2}
\nabla_X\psi_{n-1} + i X^{0,1}\cdot\psi_n = 0,\ \   
\nabla_{X^{1,0}}\psi_n + i X^{1,0}\cdot\psi_{n-1} = 0, 
\ \ \textrm{ for any } X,
\end{equation}
for the pair $(\psi_{n-1}=\bar\omega\otimes\psi_0,\psi_n= \frac{1}{2in}
\mathcal{D}\psi_{n-1})$,
and, whenever $n\neq 1$, this is finally equivalent to the K\"ahlerian 
Killing spinor equations \cite{kirchberg-killing-kahler}.
Let $\psi_0 = \left( dz_1\wedge \ldots \wedge dz^m \right)^{n/(m+1)}$ and
$f = 1 - \sum |z^i|^2$ on $\CM\mathbf{H}^m$ seen as the unit ball.
Then, Kirchberg's Ansatz \cite{kirchberg-killing-kahler} shows that, for 
any choice of multi-index 
$\alpha = (\alpha_{i_1}, \ldots, \alpha_{i_{n-1}})$, $\omega=f^{-n} 
\, dz^{\alpha_{i_1}}\wedge \ldots\wedge dz^{\alpha_{i_{n-1}}}$
provides a K\"ahlerian Killing spinor, as does $\omega =
f^{-n}\, \iota_{\sum (z^j\partial_j)} 
(dz^{\beta_{i_1}}\wedge \ldots\wedge dz^{\beta_{i_{n}}})$
for any choice of multi-index $\beta = (\beta_{i_1}, \ldots,\beta_{i_n})$.

In complex dimension $m=2$ ($n=1$), Kirchberg's Ansatz above does not 
provide
{\sl a priori} K\"ahlerian Killing spinors (the second set of equations
above is not {\sl in general} equivalent to the K\"ahlerian Killing 
condition) but it is easily checked that the spinor fields given a few lines
above are {\sl indeed} K\"ahlerian Killing spinors. \end{proof}

\begin{remk} Proposition \ref{pro:exiskki} is part of a more general 
phenomenon: for each $r$ in $\{1,\ldots,m\}$, there exists on the complex 
hyperbolic space a well chosen $\spic$--structure, built from a root of 
the canonical bundle, endowed with distinguished spinor fields living 
in $\Sigma^c_{r-1}\oplus \Sigma_r^c$. This remark can be easily 
substantiated by using Kirchberg's Ansatz.\end{remk}

\begin{remk} Our special spinors can also be obtained by projecting (in the
sense of \cite{moroianu-dirac,moroianu-complex})
the parallel spinors of $\CM^{2n+2}$, tensored with an adequate root
of the complex volume form, over the complex hyperbolic space. This idea
leads to the proofs developed in the next sections. \end{remk}

\medskip

We now come back to our general manifold $(M,g,J)$ of complex dimension
$m=2n$ equipped with its distinguished $\spic$--structure.
If $\nabla$ and $\dn$ denote the Levi-Civita connection and Dirac operator 
on $\Sigma^c$, we also define a modified connection $\nuv$ and a modified 
Dirac operator $\dl$ acting on a spinor $\psi = \sum \psi_q$ as
\begin{align*} \nuv_X \psi & \ = \ \nabla_X\psi \ + \ i\, X^{1,0}\bu 
\psi_{n-1} \ + \ i\, X^{0,1}\bu \psi_n \ , \\ 
\dl\psi & \ = \ \dn\psi \ - \ i\, m\, \sum_{q-n\ \mathrm{even}} \psi_q \ - 
\ i\, (m+2)\, \sum_{q-n\ \mathrm{odd}} \psi_q\ .\end{align*}
It is crucial to notice at this point that the modified Dirac operator
$\dl$ {\sl is not} the Dirac operator naturally issued from the modified
connection $\nuv$ (hence the difference in the notation). The Dirac operator
issued from $\nuv$ is generally not coercive on $L^2$ and is then useless,
whereas the modified Dirac operator $\dl$ is coercive, at least in our
situation. This discrepancy plays a major role in the arguments 
below.

The main tool in this section is the Weitzenb\"ock fomula for spinors 
proven in \cite[section 3]{mh-hyper}. We now rewrite it in the $\spic$
context, with the above choice of coefficients adapted to 
our needs. We also correct two typing mistakes in \cite{mh-hyper}.
 
\begin{lemm}\label{lem:blw}
For any spinor $\psi$, let $\alpha^{\psi} (X) = 
\la  \nuv_X\psi + X\bu\dl\psi, \psi\ra$. Then
\begin{align*}
- \div\, \alpha^{\psi} \ = \ & | \nuv\psi |^2 - |\dl\psi|^2 \ + \ 2i \,
\sum_{q=0}^{m} \, (-1)^{q-n}\, 
\la ( \dl\psi )_{q} , \psi_q \ra \\
& + \sum_{q=0}^{m} \,\la \,\left( \frac{1}{4}\Scal + m(m+2) - 
2(m-q)(u_q)^2 - 2q (v_q)^2 - \frac{1}{2}F\bu \right)\psi_q , \psi_q \ra .
\end{align*}
\end{lemm}

\begin{lemm} The bundle $\Sigma_{n-1}^c\oplus\Sigma_n^c$ is trivialized by
$\nuv$-parallel spinors.
\end{lemm}

\begin{proof} The first step is to prove that there exists a full set of
$\dl$-harmonic spinors on $M$, asymptotic to the model spinors on 
$\CM\mathbf{H}^m$ described in the previous section; once these spinors 
have been obtained, the proof will go along the arguments of 
\cite{mh-hyper}. 
The key point in this first step is to show that the zeroth order terms in 
the Weitzenb\"ock formula in Lemma 3.1 are always nonnegative if the 
curvature assumptions of the theorem are satisfied. Taking into account 
$F= -2i\omega$, and letting $\Scal = -4m(m+1)\kappa$ and $p = q - n$, 
the zeroth order term acting on sections of $\Sigma^c_q$ is
\[ -m(m + 1) \kappa  + m(m+2) - 2 \left(n - p \right) 
u_{n+p}^2  - 2\left(n + p\right) v_{n+p}^2 - 2p .\]
A case by case check yields that this is always
equal to $m(m+1)(1-\kappa)$, hence nonnegative.
It then implies that the modified Dirac operator $\dl$ is coercive in $L^2$,
so that we may find $\dl$-harmonic spinors for the $\spic$--structure on 
$M$ 
which are $L^2$ perturbations of the model K\"ahlerian Killing spinors at 
infinity. 

As a second step, we apply the Weitzenb\"ock formula  to any such spinor 
field $\Psi$ and get 
\[
\lim_{r\rightarrow\infty} \int_{S_r} \alpha^{\Psi} (\nu_r) 
\ \ \geq \ \ \int_M | \nuv \Psi |^2. 
\]
Arguing as in \cite{mh-hyper}, our asymptotic conditions yield that the 
limit is zero, so that $\Psi$ is $\nuv$-parallel. 

We now write $\Psi =\sum \psi_q$ with $\psi_q$ in $\Sigma_q^c$. 
As $\nuv$ respects the
splitting of $\Sigma^c$ into its $\{\psi_{n-1},\psi_n\}$--components (where 
it is a modified connection) and the remaining components $\{\psi_q\ , \ 
q\neq n-1,n\}$ (where it equals the Levi-Civita connection), one gets 
that $\nuv(\psi_{n-1} + \psi_n)=0$ and that each component $\psi_q$ is 
parallel if $q\neq n-1,n$. These remaining components are zero since they 
are in $L^2$ (model K\"ahlerian Killing spinors live in 
$\Sigma^c_{n-1}\oplus\Sigma^c_n$). We finally get
that each solution $\Psi = \psi_{n-1} + \psi_n$ solves the K\"ahlerian 
Killing equation.
As model spinors trivialize on $\CM\mathbf{H}^m$, the solutions
trivialize $\Sigma^c_{n-1}\oplus\Sigma^c_n$ as well.
\end{proof}

\smallskip

\begin{lemm} $(M,g,J)$ is K\"ahler-Einstein, with Ricci form
$\rho = -2(m+1)\omega$.
\end{lemm}

\begin{proof} Let $\Psi = \psi_{n-1}+\psi_n$ be any of our special spinors 
and $X$ a tangent vector. We compute in two different ways:
\begin{equation}\begin{split}
m X^{1,0}\bu \Psi + (m+2)X^{0,1}\bu\Psi & = \sum_{i=1}^{2m} e_{i}\bu
\left( - \nabla_{e_{i}}\nabla_{X} + \nabla_{X}\nabla_{e_{i}}\right)\Psi\\
& = -\frac{1}{2}\Ric(X)\bu\Psi - i \sum_{i=1}^{2m} \omega(e_{i},X)
e_{i}\bu\Psi\\
& = -\frac{1}{2}\Ric(X)\bu\Psi + i JX \bu\Psi .
\end{split}\end{equation}
Letting $Z = (\Ric(X) +2(m+1)X)^{1,0}$, then $Z\bu\psi_{n-1} =0$ and 
$\overline{Z}\bu\psi_{n} =0$. Since, at any point, either $\psi_{n-1}$ 
or $\psi_n$ is non zero, this forces $Z=0$. \end{proof}

\begin{remk} A direct consequence is a proof of the Theorem in complex
dimension $m=2$. The preceding arguments provide a full set 
of spinors trivializing $\Sigma_0^c\oplus \Sigma^c_1$ and the metric is 
K\"ahler-Einstein. There is on $M$ another (elementary) $\spic$--structure,
associated with the choice of determinant bundle $L'=K_M$, the canonical 
bundle of $M$. The associated spinor bundle $\Sigma'{}^c$ has a trivial 
parallel section: the volume form, that generates $\Sigma'{}^c_2$.
One may now argue algebraically and show that the
metric has constant negative holomorphic sectional curvature. 

This method does unfortunately not work in higher dimensions: when one 
seeks 
distinguished spinor fields living in $\Sigma^c_{r-1}\oplus\Sigma^c_r$ for 
arbitrary $r$, coercivity may fail in the above Weitzenb\"ock formula,
even if one looks for different $\spic$--structures.\end{remk}

\medskip

\section{The circle bundle and its cone}

We are now ready for the remaining parts of the proof : our aim is to 
show that our manifold has constant (negative) holomorphic sectional
curvature. As in \cite{mh-hyper}, we use the construction devised by 
Ch.~B\"ar \cite{baer-holonomy} and A.~Moroianu 
\cite{moroianu-dirac,moroianu-complex}.

\begin{lemm} The total space of the auxiliary bundle $P$ has a Lorentz 
Einstein metric, invariant by the $S^1$-action. The $\spic$--structure on 
$M$ lifts to a $\spin$--structure on $P$.
\end{lemm}

\begin{proof} We define the metric $g_P$ on $P$ such that horizontal and
vertical spaces given by the Levi-Civita connection are orthogonal, the 
bundle projection $\pi$ becomes a metric-preserving submersion and the 
vertical vector field $V$ (induced from the $S^1$-action) has norm $-1$. 
From Remark 9.78 in \cite{besse} it is an Einstein metric. We now follow 
the argument in \cite{moroianu-complex}. We extend the $\spic$
frame bundle $P_{\spic}M$ to a $\spic(n,1)$ principal bundle $Q$ by letting 
$Q = P_{\spic}M \times_{\spic(n)}\spic(n,1)$. The bundle $\pi^*Q$ is a
$\spic$ frame bundle defining a $\spic$--structure on $P$. Moreover, the
auxiliary bundle of this structure is the pull-backed bundle $\pi^*P$ 
which is trivial. Any choice of a global section then induces a $\spin$ 
structure on $P$.\end{proof}

\smallskip

Let $C = \RM^*_+\times P$ be the cone over $P$, endowed with the metric
$g_C = -ds^2 + s^2 g_P$ and complex structure $J_C$ given by the lift of 
the complex structure $J$ on $C$ and by the relation $J_C V = (-1)^{n} s 
\frac{\partial}{\partial s}$. Its K\"ahler form will be denoted by 
$\omega_C$.

\begin{lemm} $(C,g_C,J_C)$ is a Ricci-flat K\"ahler spin manifold of 
signature $(2m,2)$. It bears a $C_{2n+1}^n$-dimensional space of parallel 
spinors that trivialize the eigen-sub\-bundle $\mathcal{S}_n=
\ker(\omega_C\cdot + i)$ of its spinor bundle.
\end{lemm}

\begin{proof} The first claim follows from straightforward computations,
using the O'Neill and Gauss-Codazzi-Mainardi formulas.
We must now understand how our K\"ahlerian Killing spinors transform when 
pulled back to $P$, then to $C$. 

We apply the computations of \cite{moroianu-dirac,moroianu-complex} to our 
case, {\sl cum grano salis} due to the Lorentz signature. The pulled-back 
bundle $\pi^*\Sigma^c$
is made into a Dirac bundle over $P$ if one fixes $V\bu\varphi = (-1)^n\,
\overline{\varphi}$ (conjugation of spinors on an even-dimensional 
manifold). The 
connection $\bar\nabla$ induced by the Levi-Civita connection on $P$ and 
the pull-backed connection on the auxiliary bundle of the 
$\spic$--structure 
acts then as follows: for any pull-backed spinor field $\pi^*\varphi$ and
any horizontal lift $X^H$ on $P$ of a vector $X$ on $M$,
\[ \bar\nabla_X^H\pi^*\varphi  = \pi^*\left(\nabla_X\varphi + 
\frac{(-1)^n}{2}\,\pi_*\nabla_{X^H} V \bu\overline{\varphi}\right).\]
As the metric $g_P$ is defined from the natural connection on $P$ whose
curvature is $F = -2i\omega$, the usual $6$-term formula yields
$\pi_*\nabla_{X^H} V = - JX$ and finally
\[ \bar\nabla_X^H\pi^*\varphi  = \pi^*\left(\nabla_X\varphi + 
\frac{(-1)^n}{2}\, JX \bu\overline{\varphi}\right).\]
Similarly, one computes
\[ \bar\nabla_V\pi^*\varphi = - \frac{1}{2}\,
\pi^*\left(\omega\bu\varphi\right) .\]
Using the canonical global section on $\pi^*P$ over $P$, we can now turn 
any spinor field for the $\spic$--structure into a section for the $\spin$
structure whose existence has been already remarked. 
Hence, if $\varphi$ is any spinor field,
it may be seen as associated to the $\spic$--structure (and is acted upon 
by the connection $\bar\nabla$ induced from the Levi-Civita connection
tensorized with the pull-back connection) or to the $\spin$--structure (and 
is then acted upon by the Levi-Civita connection $\nabla^P$ which is better 
seen as Levi-Civita tensorized by the trivial connection). But the trivial
connection on $\pi^*P$ differs from the pulled back connection by a factor
which is exactly the connection $1$-form $\alpha$ of $P\rightarrow M$ and 
this shows that $\nabla^P_{X^H}\varphi = \bar\nabla_{X^H}\varphi$ and 
$\nabla^P_V \varphi = \bar\nabla_V\varphi - \frac{i}{2} \varphi$.

We now apply the preceding remarks to our set of K\"ahlerian Killing spinors
$\Psi$ trivializing $\Sigma^c_{n-1}\oplus\Sigma^c_n$ on $M$. From Formula
(\ref{eq:kki2}), we get that the pull back spinors $\Psi^P$ solve
\begin{equation}\label{eq:ki}
\nabla^P_{W}\Psi^P + \frac{i}{2}W\bu\Psi^P\ =\ 0\ \ \ \ \forall\, W \in TP.
\end{equation}
And it is now straightforward to show that these can be extended as
parallel spinors on the cone $C$ living in the desired subbundle.
\end{proof}

\begin{lemm} $C$ is flat and $M$ has constant negative holomorphic 
sectional curvature.
\end{lemm}

\begin{proof} 
The previous Lemma provides parallel spinors trivializing the component 
${\mathcal S}_n$ (with respect to the action of the K\"ahler form of $N$) 
of the spin bundle $\mathcal{S}$ of $C$. The Weyl tensor $W^C$ acts
trivially on ${\mathcal S}_n$, {\sl i.e.} for each $X$, $Y$ in $TC$, 
\begin{equation} W_{X,Y}^C \cdot \psi = 0 .\end{equation}
Mimicking the representation-theoretic argument in \cite{mh-hyper} 
shows that the only $2$-form
in $\Lambda^{1,1}_0C$ whose action on ${\mathcal S}_n$ is zero is the
zero form itself. 
This proves $C$ is a flat manifold; applying O'Neill formulas yields
that $M$ has constant negative holomorphic sectional curvature 
\cite[Corollary 9.29]{besse}. \end{proof}

\begin{lemm} $(M,g,J)$ is biholormophically isometric to $\CM\mathbf{H}^m$.
\end{lemm}

\begin{proof} We only need to check that $M$ is simply-connected. Assume 
the contrary : let $M_0 = \CM\mathbf{H}^m$ be the universal 
covering of $M$. If the covering is not trivial, $M_0$ has at 
least two ends since the unique end of $M$ is simply-connected. This is
of course a contradiction. \end{proof}

\bigskip

\begin{small}
{\flushleft\it Acknowledgement}. We thank A.~Moroianu and T.R.~Ramadas for 
some useful discussions about this paper.
\end{small}

\newpage

\phantom{...}

\bigskip

\bibliographystyle{amsplain}
\providecommand{\bysame}{\leavevmode\hbox to3em{\hrulefill}\thinspace}

\bigskip

\end{document}